\documentclass{icmart}

\usepackage[all]{xy}

\contact[bertrand.toen@um2.fr]{Universit\'e de Montpellier 2, Case courrier 051, B\^at 9, Place Eug\`ene 
Bataillon, 
Montpellier Cedex 5, France}

\newcommand {\Map} {\mathbf{Map}}
\newcommand {\Parf} {\mathbf{Perf}}

\newcommand {\OO} {\mathcal{O}}

\newcommand {\Spec} {\mathbf{Spec}}

\newcommand  {\cdga}     {\mathbf{cdga}}
\newcommand  {\art}     {\mathbf{dg-art}^{*}}

\newcommand  {\dSt}   {\mathbf{dSt}}
\newcommand  {\dSch}   {\mathbf{dSch}}
\newcommand  {\dArSt}   {\mathbf{dArSt}}

\newcommand  {\Top}   {\mathbb{S}}

\newcommand{\s}{\infty}

\newcommand  {\Pol}   {\mathcal{P}ol}

\newtheorem{thm}{Theorem}[section]
\newtheorem{princ}{Principle}[section]
\newtheorem{cor}[thm]{Corollary}

\newtheorem{conj}[thm]{Conjecture}

\newtheorem{df}[thm]{Definition}

\title[Derived Algebraic Geometry and Deformation Quantization]{Derived Algebraic Geometry and Deformation 
Quantization}

\author[Bertrand To\"en]{Bertrand To\"en}

\begin{document}

\begin{abstract}
This is a report on recent progress concerning the interactions between derived algebraic geometry and 
deformation quantization. 
We present the notion of derived algebraic stacks, of shifted symplectic and Poisson structures, as well as the 
construction of 
deformation quantization of shifted Poisson structures. As an application we propose a general construction of 
the quantization
of the moduli space of $G$-bundles on an oriented space of arbitrary dimension.
\end{abstract}

\maketitle

\section{Introduction}

Quantization is an extremely vast subject, particularly because it has 
a long-standing physical origin and history. 
Even from the more restrictive point of view of a pure mathematician,
quantization possesses many facets and connects with a wide variety of modern mathematical domains.
This variety of interactions explains the numerous mathematical incarnations that the expression 
\emph{quantization} finds in the existing literature, though a common denominator seems 
to be a \emph{perturbation of a commutative structure into a non-commutative structure}. 
For a commutative object $X$ (typically a commutative algebra or
a manifold), a quantization is most often realized as a family $X_{\hbar}$, of objects
depending on a parameter $\hbar$, which recovers $X$ when $\hbar=0$ and which is non-commutative
for general values of $\hbar$. The existence of the family $X_{\hbar}$ 
is in most cases related to the existence of certain additional geometric structures, such as symplectic or 
Poisson structures. 

The purpose of this manuscript is to present a new approach to quantization, or more specifically
to the construction and the study of interesting non-commutative deformations of commutative objects of
geometrico-algebraic origins. This new approach is based on the \emph{derived algebraic geometry}, 
a version of algebraic geometry that has emerged in the last decade (see \cite{seatt,dagems}), and which itself 
consists of
a \emph{homotopical perturbation of algebraic geometry}. Derived algebraic geometry 
not only leads to a unified geometric interpretation of most of the already existing
quantum objects (e.g. it treats the quantum group and deformation quantization of a Poisson manifold on 
an equal footing), 
but also opens up a whole new world
of quantum objects which, as far as we know, have not been identified in the past even though they seem 
to appear naturally in algebraic geometry, algebraic topology, or representation theory. \\

\textbf{Acknowledgements:} The content of this manuscript represents works still in progress in 
collaboration with T. Pantev, M. Vaqui\'e and G. Vezzosi. The author is thankful to the three of them
for numerous conversations that
have helped a lot in the writting of this manuscript. An important part of the present text 
has been prepared during a, snowy but enjoyable, visit at Yale Mathematics department. A particular thanks
to S. Goncharov and M. Kapranov for many enlightening conversations on related subjects. \\

\textbf{Convention:} All varieties, algebras, schemes, 
stacks, algebraic groups etc \dots will be over a ground field $k$ of characteristic 
zero. 

\section{Quantization as deformed categories: three motivating examples}

In this first section we briefly recall three well known important examples of "quantization in action" in 
different
domains: 
quantum groups, skein algebras and Donaldson-Thomas invariants. We identify the natural moduli spaces 
behind each of these examples 
and  explain how they all can be considered from the unified point of view of deformations of categories
and monoidal categories of sheaves. \\

\textbf{Quantum groups.} Probably the most famous and most fundamental of quantum objects are quantum groups. For 
an algebraic group $G$, with lie algebra $\frak{g}$, and a choice of a $G$-invariant element $p \in 
Sym^{2}(\frak{g})$,
Drinfeld constructs a quantum group (see \cite{dr}). Algebraically the quantum group is a deformation 
of the Hopf algebra $A=\mathcal{O}(G)$ of functions on $G$, into a non-commutative Hopf algebra
$A_{\hbar}$.  

\textbf{Skein algebras.} Skein algebras appear in low dimensional topology (see \cite{tu}). They are
associated with a given Riemann surface $\Sigma$, and are explicitly defined in terms of generators and 
relations. 
The generators are given by simple curves traced on the surface $\Sigma$, and the relations are given by the
so-called skein relations, which possess natural deformations by a parameter $q=e^{2i\pi\hbar}$. 
The skein algebra associated with $\Sigma$, $K_{\hbar}(\Sigma)$, is a non-commutative deformation
of the ring of functions on the character variety 
of $\Sigma$ for the group $Sl_{2}$ (i.e. the affine algebraic variety whose points describe
$Sl_2$-representations of the fundamental group $\pi_1(\Sigma)$). 

\textbf{Donaldson-Thomas invariants.} For $X$ a Calabi-Yau algebraic variety of dimension $3$, 
we denote by $\mathcal{M}_{X}$ the moduli space of stable vector bundles with fixed numerical invariants. 
It is a singular variety in general but with 
a very specific local structure. Indeed, it is known that locally around each point, $X$
embeds into a smooth ambient variety $Z$ as the critical points of a function
$f : Z \rightarrow \mathbb{A}^{1}$. Each of these locally defined functions $f$ define
a (perverse) sheaf $\nu_{f}$ of vanishing cycles on $X$, which under an
orientability assumption glue to a globally defined 
perverse sheaf $\mathcal{E}$ on $X$ (see \cite{bbbj,bbdjs,bu} for more on the subject). The sheaf 
$\mathcal{E}$ is a quantization of the space $X$, in the sense that 
it can be seen to be a deformation of the line bundle of \emph{virtual half forms on $X$} 
(we refer here to the next section for more about virtual structures). This
deformation is again a non-commutative deformation, but this time in a dramatic 
way as the multiplicative structure itself is lost and  $\mathcal{E}$ only exists
as a sheaf of (complexes of) vector spaces on $X$. \\

Despite their different origins and  
differences in appearance these three
examples of quantization can be considered in a striking unified way: they all are 
deformations of categories of sheaves on natural moduli spaces, where the categories 
eventually come equipped with monoidal structures. More is true, these deformations are all 
induced by the same type of structures on the corresponding moduli spaces, at least when 
they are appropriately viewed as derived algebraic stacks as we will explain later.
The moduli spaces related to these three examples are easy to guess: they are respectively the
moduli space $Bun_{G}(*)$ of $G$-bundles on a point $*$, the moduli space $Bun_{Sl_2}(\Sigma)$ 
of $Sl_2$-bundles on the surface $\Sigma$, and the moduli space $\mathcal{M}_{X}$ of 
algebraic $Gl_n$-bundles on $X$, also denoted $Bun_{Gl_n}(X)$. In the case of the quantum group,
there is no non-trivial $G$-bundles on a point, but the trivial $G$-bundle possesses many 
automorphisms. The moduli space $Bun_{G}(*)$ is thus trivial from the point of view of algebraic
varieties but can be realized as a non-trivial algebraic stack $BG$. Quasi-coherent sheaves on 
$BG$ are nothing else than linear representations of $G$, $QCoh(BG)=Rep(G)$. The quantum
group $A_{\hbar}$ can then be realized as a deformation of $QCoh(BG)$ considered as
a braided monoidal category. For the case of skein algebras, 
we already mentioned that $K_{\hbar=0}(\Sigma)$ is the ring of functions on the moduli space
$\chi(\Sigma)=Bun_{Sl_2}(\Sigma)$. The moduli space $Bun_{Sl_2}(\Sigma)$ is an affine algebraic
variety and thus its category of quasi-coherent sheaves is equivalent to modules over 
its ring of functions, $QCoh(Bun_{Sl_2}(\Sigma)) = K_{\hbar=0}(\Sigma)-Mod$. The deformation
$K_{\hbar}(\Sigma)$ can thus be realized as a deformation of the category 
$QCoh(Bun_{Sl_2}(\Sigma))$, simply considered as a linear category. Finally, 
the perverse sheaf $\mathcal{E}$ on $\mathcal{M}_{X}$ can itself be considered as a deformation of 
a natural object $\omega_{X}^{1/2,virt}$, of \emph{virtual
half top forms on $X$}, which is almost a quasi-coherent sheaf on $X$ (it is
a complex of such). The quantized object $\mathcal{E}$ is thus not
a deformation of $QCoh(X)$, but is rather a deformation of one of its objects. 

To summarize, all of the three examples discussed above have an interpretation in terms of deformations
of categories, possibly with monoidal structures, of quasi-coherent sheaves on certain moduli spaces. 
Monoidal categories can be organized in a hierarchy, corresponding to the degree of symmetry 
imposed on the monoidal structure. For instance, a monoidal category can come equipped with a braiding, or a 
symmetry constraint. Monoidal categories will be referred to as $1$-fold monoidal categories, 
braided monoidal categories as $2$-fold monoidal categories, and symmetric monoidal
categories as $\infty$-fold monoidal categories. We will moreover see that when categories are 
replaced by $\infty$-categories there is a notion of $n$-fold monoidal 
$\infty$-categories for $2 < n < \infty$
(also called $E_{n}$-monoidal $\infty$-categories), 
interpolating between braided and symmetric monoidal categories. When $n=0$, a $0$-fold category
can be defined to simply be a category,  a $(-1)$-monoidal category can be declared to be 
an object in a category, and a $(-2)$-monoidal category can be defined as 
an endomorphism of an object in a category. This hierarchy is rather standard in the setting of higher category
theory in which a monoidal category is often considered as a 2-category with a unique object, and
a braided monoidal categories as a 3-category with unique object and unique $1$-morphism 
(see for instance \cite[\S V-25]{si1}).

In our examples above, quantum groups are deformations of $QCoh(Bun_{G}(*))$ as 
a 2-fold monoidal category. Skein algebras are deformations of 
$QCoh(Bun_{Sl_2}(\Sigma))$ as $0$-fold monoidal categories. Finally, the perverse sheaf
$\mathcal{E}$ is a deformation of $QCoh(Bun_{Gl_n}(X))$ considered
as $(-1)$-fold monoidal category. The purpose of the present paper is to explain that this is only
a very small part of a bigger coherent picture, which we present here as a key principle.

\begin{princ}\label{princ} 
For any oriented manifold of dimension $d$ (understood either in the topological
or in the algebraic sense), and any reductive group $G$, the moduli space of $G$-bundles
on $X$,
$Bun_{G}(X)$, possesses a quantization which is a deformation of $QCoh(Bun_{G}(X))$
considered as an $(2-d)$-fold monoidal $\infty$-category. 
\end{princ}

We will see how this principle can becomes a theorem, after a suitable interpretation
of $Bun_{G}(X)$, $QCoh(Bun_{G}(X))$, and a suitable understanding of $(2-d)$-fold monoidal structures. 
We will also see how this principle follows from the general framework of symplectic and poisson structures
in derived algebraic geometry, and a general quantization procedure. 

\section{Moduli spaces as derived stacks}

The concept of  an algebraic variety is not 
quite enough to encompass all the aspects of the moduli problems appearing in algebraic
geometry. Starting in the 50' and continuing until this day, several successive 
generalizations of algebraic varieties were introduced in order to understand 
more and more refined aspects of moduli spaces. As a first step nilpotent functions have been 
allowed as it is well known that many interesting moduli 
spaces are non-reduced and must be considered as schemes instead of algebraic varieties. 
Secondly algebraic stacks have been introduced in order to take into account the fact that in most examples, 
moduli spaces classify objects only up to isomorphisms, and in many situations non-trivial automorphisms
prevent the existence of any reasonable  moduli spaces. Unfortunately algebraic stacks are still not
enough to capture all aspects of moduli problems, as even though they see non-trivial automorphisms 
the so-called \emph{higher structures} remain invisible. We will explain in this section what the
higher structures are and how the notion of \emph{derived algebraic stacks} is needed
in order to incorporate them as part of the refined moduli space. 

\subsection{Higher structures I: higher stacks} 
A first type of 
higher structure concerns \emph{higher homotopies}, which appear naturally each time 
objects are classified not only up to isomorphism but up to a weaker notion of
\emph{equivalences}. A typical example is the extension of the moduli space of
vector bundles on a given smooth and projective algebraic variety $X$, by also allowing
complexes of vector bundles, now considered up to quasi-isomorphism\footnote{This appears typically 
in Donaldson-Thomas theory for which moduli spaces of objects in the bounded coherent 
derived category $D^{b}_{coh}(X)$ must be considered.}. 
The moduli space of 
vector bundles on $X$ can be realized as an algebraic stack, but the moduli 
of complexes of vector bundles taken up to quasi-isomorphism can not be represented by  
an algebraic stack in the sense of \cite{ar}. The reason for this is the existence of higher homotopies
between maps between complexes, which is reflected in the fact that a complex $E$ on $X$ can have 
non-trivial negative self extension groups $Ext^{-i}(E,E)$. The vector spaces
$Ext^{-i}(E,E)$ for $i>0$ are higher analogues of automorphism groups of vector bundles and their
existence prevent the representability by an algebraic stack for the exact same reason 
that the existence of non-trivial automorphisms of vector bundles prevent the representability
of the moduli problem of vector bundles on $X$ by a scheme. 
In his manuscript "Pursuing stacks", A. Grothendieck brought forward the idea of \emph{higher stack}, which 
is an extension of the notion of stacks of groupoids usually considered in moduli theory
to a higher categorical or higher homotopical setting. This idea has been made concrete  in \cite{si2}
by the introduction of a notion of \emph{algebraic $n$-stacks} (see also \cite{seatt}). These algebraic $n$-stacks
behave in a very similar way to algebraic stacks, and most of the standard notions and 
techniques of algebraic geometry 
remain valid in this new setting (they have derived categories, cohomology, tangent 
spaces, dimensions \dots). Fundamental examples of algebraic $n$-stacks include
the Eilenberg-McLane stacks of the form $K(A,n)$, for $A$ a commutative algebraic group, which 
are higher analogues of classifying stacks $BG$. Another important example for us are the so-called linear 
stacks: for a scheme $X$ and a complex of vector bundles $E_*$ on $X$ concentrated in degrees
$[-n,0]$, there is a linear stack $\mathbb{V}(E_*) \longrightarrow X$, which is
a generalization of the total space of a vector bundle. 
Finally, for $X$ a smooth and projective variety, there is an algebraic $n$-stack 
of complexes of vector bundles on $X$, which also possesses many possible 
non-commutative generalizations (see \cite{seatt,dagems} for more on the subject). \vspace{5mm}

\subsection{Higher structures II: derived algebraic stacks} A second type of higher structure
attached to moduli problems is called the \emph{derived structure}. These derived structures 
are somehow dual to the
higher homotopies we have just mentioned and exist even in absence of any stacky 
phenomenon (i.e. even when there are no non-trivial automorphisms). They have been introduced 
through the eye of deformation theory and originally were only considered at the formal
level of moduli spaces. The \emph{derived deformation theory}, also referred to as DDT, 
is a collection of ideas going back to the 80's, stipulating that moduli spaces, formally around
a given fixed point, can be described in terms of Mauer-Cartan elements in a suitable dg-Lie algebra
associated to this point. The most famous example is the deformation theory of 
a given smooth projective variety $X$, for which the natural dg-Lie algebra is $C^*(X,\mathbb{T}_{X})$, the
cochain complex computing the cohomology of the tangent sheaf, endowed with its dg-Lie structure
coming from the bracket of vector fields. This example is not special, 
and in fact all possible moduli problems
come with natural dg-Lie algebras describing their formal completions. \\

A striking consequence of the DDT is the existence of \emph{virtual sheaves} on moduli spaces. Indeed, 
according to the DDT, for a point $x\in \mathcal{M}$ in some moduli space $\mathcal{M}$, we can find a dg-Lie
algebra $\frak{g}_{x}$ controlling the formal local ring of $\mathcal{M}$ at $x$. There 
is moreover an explicit formula reconstructing formal functions at $x$:
$\widehat{\OO}_{\mathcal{M},x} \simeq H^{0}(C^{*}(\frak{g}_{x}))\simeq H^{0}(\frak{g}_{x},k),$
where $C^{*}(\frak{g}_{x})=\widehat{Sym}_{k}(\frak{g}_{x}^*[-1])$ is the (completed)
Chevalley complex of $\frak{g}_{x}$, which also computes the cohomology of 
$k$ considered as a the trivial $\frak{g}_{x}$-module. An important observation is that the 
Chevalley complex $C^*(\frak{g}_{x})$ is a commutative dg-algebra which can have
non trivial cohomology in non-positive degrees. These cohomology groups, $H^{i}(C^*(\frak{g}_{x}))$
for $i<0$, provide non-trivial coherent sheaves over the formal neighborhood of $x$, which 
are by definition the \emph{derived structures of $\mathcal{M}$ around $x$}. These local coherent 
sheaves are quite important, as they control for instance the smoothness defect of the
moduli space $\mathcal{M}$, and lead to the so-called virtual fundamental class (see \cite{ko1}). 
Incorporating these higher structures as an intrinsic part of the moduli space itself
has lead to the theory of \emph{derived algebraic geometry}, and to introduction of 
\emph{derived schemes} and \emph{derived algebraic ($n-$)stacks} as the correct
geometrico-algebraic notion to fully represent moduli problems in algebraic geometry.

\subsection{Derived schemes and derived algebraic stacks} The foundations of the theory of derived algebraic geometry can be found in \cite{hagI,hagII} and \cite{lu1}. 
We will not give precise definitions here, as the details easily become technical, 
and will rather concentrate on some basic definitions and basic facts. 

\subsubsection{Derived schemes} As objects derived schemes are rather
easy to define and understand. We display below one possible definition of derived schemes (specific
to the characteristic zero case, recall that everything is over a base field $k$ of zero characteristic).

\begin{df}\label{dsch}
A \emph{derived scheme} (over the field $k$) consists of a pair $(X,\OO_{X})$, where
$X$ is a topological space and $\OO_{X}$ is a sheaf of commutative differential graded $k$-algebras 
on $X$, satisfying the following conditions.
\begin{itemize}

\item The sheaves $H^{i}(\OO_{X})$ vanish for $i>0$. 
\item The ringed space $(X,H^{0}(\OO_{X}))$ is a $k$-scheme.
\item For all $i$, $H^{i}(\OO_{X})$, considered as a sheaf of $H^{0}(\OO_{X})$-modules, 
is a quasi-coherent sheaf. 

\end{itemize}
\end{df}

The above definition makes derived schemes look as rather simple objects, but things get more sophisticated 
when morphisms between derived schemes are introduced. The sheaf of dg-algebras $\OO_{X}$ 
must only be considered up to quasi-isomorphisms, and quasi-isomorphic derived schemes have to be 
considered as equivalent. Therefore, there is a \emph{derived category of derived schemes}, which 
is non-linear analogue of the derived category of a ring, and for which quasi-isomorphic
sheaves of dg-algebras define the same derived scheme. There are two possible constructions 
of the category of derived schemes, a first one relies on model category structures and 
includes the quasi-isomorphisms as the weak equivalences of a certain model category.
A more modern approach, powerful in practice but 
more demanding in terms of foundations, is to use $\infty$-categories, and to define
the category of derived schemes directly as an $\infty$-category (see for instance 
\cite{si1}). Concretely this means that 
morphisms between two given derived schemes do not form a set anymore but are topological spaces, 
or a simplicial
sets. This allows to consider homotopies between morphisms of derived schemes, and thus to define
equivalences between derived schemes as morphisms having inverses up to homotopy. The $\s$-category 
$\dSch$ of derived schemes is then defined so that quasi-isomorphisms become homotopy equivalences
in $\dSch$. We refer to \cite[\S 2.1]{dagems} for more details on these two approaches, 
and we will consider $\dSch$ as an $\s$-category in what follows. 

There have been a certain number of works on the notion of derived schemes, 
making many of the basic aspects of scheme theory available in the derived setting. Derived schemes
behave in a very similar fashion to schemes, they have a notion of (quasi-coherent) sheaves,
cohomology, smooth, flat and \'etale maps etc \dots. Special among  the derived schemes are 
the affine derived schemes, which are completely characterized by their functions, which themselves
form a non-positively graded cdga. There is a $\Spec$ construction, sending a 
cdga $A$ to an affine derived scheme $\Spec\, A$, whose underlying space is $\Spec\, H^{0}(A)$ and
whose structure sheaf is given by the various localizations $A[f^{-1}]$ in a very similar manner 
as for un-derived schemes. The $\Spec$ functor produces a full embedding of the (opposite) $\s$-category 
of cdga into $\dSch$. Here the $\s$-category of cdga can be presented concretely as the category 
whose objects are quasi-free cdga together with the standard simplicial sets of morphisms 
$Map(A,B)$\footnote{Whose set of $n$-simplicies is $Hom(A,B\otimes_{k} DR(\Delta^n))$, where 
$DR(\Delta^n)$ is the algebraic de Rham complex of the algebraic $n$-dimensional simplex $\{\sum x_i = 1 \} 
\subset 
\mathbb{A}^{n+1}$.}. 
A general derived scheme $X$ is locally equivalent to $\Spec\, A$ for some cdga $A$, and many of the
notions defined for cdga can be extended to arbitrary derived scheme by sheafification. 
This is for instance the case for the notions of smooth, flat and \'etale
maps, as well as for the notion of cotangent complexes of derived schemes, etc. \vspace{5mm}

\subsubsection{Derived algebraic stacks}
The reader should have already guessed that derived schemes are not quite enough for our purpose and that 
we will need the notion of derived algebraic stacks (including derived higher stacks in some cases). 
These are defined in a similar fashion as  algebraic stacks and higher algebraic stacks (see 
\cite{hagII,dagems} for details). In a nutshell a derived algebraic stack is given by a quotient 
of a derived scheme $X$ by an action of a smooth groupoid. Concretely 
a derived algebraic stack is associated to a simplicial object $X_{*}$ made of 
derived schemes satisfying some smooth Kan lifting conditions (see \cite{pr}). A typical example is the
action of an algebraic group $G$ on a derived scheme $Y$, for which the simplicial 
objects is the standard nerve of the action
$([n] \in \Delta ) \mapsto Y\times G^{n},$
where the face maps are defined by the action of $G$ on $Y$ and the multiplication on $G$, and the
degeneracies are induced by the unit in $G$. The derived algebraic stack obtained this way 
will be denoted by $[Y/G]$, and it should be noted that already some interesting derived algebraic stacks
are of this form (but these are not enough to represent all moduli problems). 

Derived algebraic stacks are good
objects to do algebraic geometry with, and many of the 
standard notions and results known for underived algebraic stacks can be extended to the derived setting.
One construction of fundamental importance for us is the (dg-)category (see \cite{ke}) 
of quasi-coherent complexes
on a given derived algebraic stack. For an affine derived scheme $X=\Spec\, A$, the 
quasi-coherent complexes over $X$ are declared to be the $A$-dg-modules. The $A$-dg-modules
form a nice $k$-linear dg-category $L(A)$, for which one explicit model consists of the dg-category
of quasi-free $A$-dg-modules. For a general derived algebraic stack $X$ the dg-category 
of quasi-coherent complexes is defined by approximating $X$ by affine derived schemes 
$$L(X):=\lim_{\Spec\, A \rightarrow X} L(A),$$
where the limit is taken inside a suitable $\s$-category of dg-categories (see \cite{to4}), or equivalently 
is understood as a homotopy limit inside the homotopy theory of dg-categories of \cite{ta}. 

Another important notion we will use is the cotangent complex. Any derived algebraic stack possesses
a canonically defined object $\mathbb{L}_{X} \in L(X)$, which is the derived version 
of the sheaf of Kalher 1-forms. When $X$ is a smooth scheme then $\mathbb{L}_{X}$ is the
vector bundle $\Omega^1_{X}$ considered as an object in the quasi-coherent derived category of $X$. 
When $X=\Spec\, A$ is an affine derived scheme, $\mathbb{L}_{A}$ is the $A$-dg-module representing
the so-called Andr\'e-Quillen homology, and can be defined as the left derived functor
of $A \mapsto \Omega^{1}_A$. For a scheme $X$, $\mathbb{L}_{X}$ coincides with Illusie's cotangent complex.
For a general derived algebraic stack $X$ the cotangent complex $\mathbb{L}_{X} \in L_{qcoh}(X)$ 
is obtained by gluing the cotangent complexes of each stage in a simplicial presentation, 
but can also be characterized by a universal property involving square zero extensions (see 
\cite[\S 3.1]{dagems}). The dual object $\mathbb{T}_{X}:=\underline{Hom}_{\OO_{X}}(\mathbb{L}_{X},\OO_{X}) 
\in L_{qcoh}(X)$ is called the tangent complex of $X$ and is a derived version of the sheaf of derivations.

To finish with general facts about derived algebraic stacks, we would like to mention a
specific class of objects which are particularly simple to describe in algebraic terms, 
and which already contains several non-trivial examples. This class consists of derived algebraic
stack of the form $X=[Y/G]$, where $Y$ is an affine derived scheme and $G$ a linear algebraic
group acting of $Y$. The derived affine scheme $Y$ is the spectrum
of a commutative dg-algebra $A$, which up to a quasi-isomorphism can be chosen 
to be a cdga inside the category $Rep(G)$ of linear representations of $G$ (and 
$A$ can even be assumed to be free as a commutative graded algebra). 
The cdga $A$, together with its strict $G$-action, can be
used in order to describe $L_{qcoh}(X)$ as well as $\mathbb{L}_{X} \in L_{qcoh}(X)$.
A model for the dg-category $L_{qcoh}(X)$ is the dg-category of 
cofibrant and fibrant $A$-dg-modules inside $Rep(G)$, where here fibrant refers to 
a model category structure on the category of complexes of representations of $G$ (and
fibrant means $\mathcal{K}$-injective complex of representations). In particular 
the homs in the dg-category $L_{qcoh}(X)$ compute $G$-equivariant
ext-groups of $A$-dg-modules. The object $\mathbb{L}_{X} \in L_{qcoh}(X)$ can be described as follows. 
The $G$-action on $A$ induces a morphism of $A$-dg-modules
$\mathbb{L}_{A} \longrightarrow \frak{g}^{\vee}\otimes_{k}A$, where 
$\frak{g}$ is the Lie algebra of $G$, which is a morphism of $A$-dg-modules inside $Rep(G)$. 
The cone of this morphism, or more precisely 
a cofibrant and fibrant model for this cone, is a model for $\mathbb{L}_{X}$ as an object in $L_{qcoh}(X)$.  

\subsection{Representability of derived mapping stacks} As for the case of derived schemes, derived 
algebraic
stacks form an $\s$-category denoted by $\dArSt$ (where "Ar" stands for "Artin"). 
This $\s$-category is itself a full sub-$\s$-category of $\dSt$, the 
$\s$-category of (possibly non-algebraic) derived stacks. The objects of $\dSt$ are $\s$-functors
$F : \cdga \longrightarrow \Top$ satisfying \'etale descent conditions, and are also called 
\emph{derived moduli problems}. The derived moduli problems can sometimes be represented by schemes, by
derived schemes, or by derived algebraic stacks, in the sense that there exists a derived algebraic 
stack $X$ together with functorial equivalences $F(A) \simeq Map_{\dArSt}(\Spec\, A,X)$. Proving that a 
given 
derived moduli problem is representable is in general not a 
trivial task, and the following theorem provides a way to 
construct new derived algebraic stacks. 

\begin{thm}{(\cite[Thm. 2.2.6.11]{hagII})}\label{t1}
Let $X$ be a smooth and proper scheme and $Y$ a derived algebraic stack which is locally of finite 
presentation 
over the base field $k$. Then, the derived moduli problem $A \mapsto Map_{\dSt}(X \times \Spec\, A,Y)$
is representable by a derived algebraic stack denoted by $\Map(X,Y)$. 
\end{thm}

One important aspect of the theorem above lies in the fact that the (co-)tangent complexes of
the derived mapping stacks $\Map(X,Y)$ are easy to compute: there is a diagram 
of derived algebraic stacks
$\xymatrix{ Y & \ar[l]_-{ev} X\times \Map(X,Y) \ar[r]^-{p} & \Map(X,Y),}$
where $ev$ is the evaluation map and $p$ is the natural projection, and we have

\begin{equation}\label{tangent}
 \mathbb{T}_{\Map(X,Y)} \simeq p_{*}ev^{*}(\mathbb{T}_{Y}) \in L_{qcoh}(\Map(X,Y)).
\end{equation}

At a given point $f\in \Map(X,Y)$, corresponding to a morphism $f : X \longrightarrow Y$, 
the formula states that the tangent complex at $f$ is $C^{*}(X,f^{*}(\mathbb{T}_{Y}))$, the cochain complex
of cohomology of $X$ with coefficients in the pull-back of $\mathbb{T}_{Y}$ by $f$. This last formula is 
moreover compatible
with the dg-Lie structures: the formal completion of $\Map(X,Y)$ corresponds, via the DDT correspondence, 
to the dg-Lie algebra $C^{*}(X,f^{*}(\mathbb{T}_{Y}))[-1]$ (here $\mathbb{T}_{Y}[-1]$ is equipped with 
its natural dg-Lie structure, see \cite{he}). This provides a nice and efficient way to understand the 
derived moduli 
space $\Map(X,Y)$ at the formal level. 

The above theorem also possesses several possible variations, which often can be reduced to the statement 
\ref{t1} itself. 
For instance $X$ can be replaced by a finite homotopy type (e.g. a compact smooth manifold), or 
by formal groupoids such as $X_{DR}$ or $X_{Dol}$ (see \cite[\S 9]{si3}). This provides existence of 
derived
moduli stacks for maps $X \longrightarrow Y$ understood within different settings (e.g. 
locally constant maps, maps endowed with flat connections or with a Higgs field etc \dots). 
The formula for the tangent complex remains correct for these variants as well, 
with a suitable definition of the functors $p_*$ and $ev^*$ involved. Pointwise this is reflected in the fact 
that 
$C^{*}(X,f^{*}(\mathbb{T}_{Y}))$ is replaced with the appropriate cohomology theory (cohomology with local 
coefficients
when $X$ is a finite homotopy type, algebraic de Rham cohomology for $X_{DR}$ etc \dots). 

The main examples of applications of theorem \ref{t1} will be for us when $Y=BG$, the classifying 
stack of $G$-bundles, in which case $\Map(X,BG)$ is by definition the derived moduli stack of $G$-bundles
on $X$. Other interesting instances of applications arise when $Y=\Parf$ is the derived stack of perfect complexes
(in which case $\Map(X,Y)$ is the derived moduli stack of perfect complexes on $X$), 
or when $Y$ is the total space of a shifted total cotangent bundle (see below).

\section{Symplectic and Poisson structures in the derived setting}

In the previous section we saw why and how moduli problems can be represented by derived schemes and 
derived algebraic
stacks. In the sequel we will be further interested in the representability of derived mapping stacks 
provided by our theorem 
\ref{t1}, as well as the formula for their tangent complexes. This formula is the key for the construction of 
symplectic 
structures on derived mapping stacks, by using cup products in cohomology in order to define pairing 
on tangent complexes. This brings us to the notions of \emph{shifted symplectic structures}, and of
\emph{shifted poisson structures}, of major importance in order to achieve the goal proposed by our 
principle \ref{princ}. 

\subsection{Algebraic de Rham theory of derived algebraic stacks}
Let $X$ be a derived algebraic stack locally of finite presentation over our ground field $k$. We have seen that $X$ possesses
a cotangent complex $\mathbb{L}_{X}$, which is a quasi-coherent complex on $X$. In our situation $\mathbb{L}_{X}$ is moreover 
a perfect $\OO_{X}$-module (see \cite[Def. 1.2.4.6]{hagII}) because of the locally finite 
presentation condition, and is thus a dualizable object 
in $L_{qcoh}(X)$. Its dual (dual here as an $\OO_{X}$-module), is the tangent complex $\mathbb{T}_{X}:=\mathbb{L}_{X}^{\vee}$.
A \emph{$p$-form on $X$} is simply defined as an element in $H^{0}(X,\wedge_{\OO_{X}}^{p}\mathbb{L}_{X})$, or equivalently 
as a homotopy class of morphisms $w : \wedge_{\OO_{X}}^{p}\mathbb{T}_{X} \longrightarrow \OO_{X}$ in $L_{qcoh}(X)$. More generally, 
if $n\in \mathbb{Z}$, 
a \emph{$p$-form of degree $n$ on $X$} is an element in $H^{n}(X,\wedge_{\OO_{X}}^{p}\mathbb{L}_{X})$, or equivalently 
a homotopy class of morphisms $w : \wedge_{\OO_{X}}^{p}\mathbb{T}_{X} \longrightarrow \OO_{X}[n]$ in 
$L_{qcoh}(X)$. For $p$ fixed, 
$p$-forms of various degrees 
form a complex of $k$-vector spaces $\mathcal{A}^{p}(X)=\Gamma(X,\wedge_{\OO_{X}}^{p}\mathbb{L}_{X})$\footnote{Here and in the sequel $\Gamma(X,-)$ stands for 
the $\s$-functor of global sections, and thus computes hyper-cohomology on $X$.}, 
whose $n$-th cohomology space is the space of $p$-forms of degree $n$.

The total complex of differential
forms on $X$ is defined as an infinite product
$$\mathcal{A}(X):=\prod_{i\geq 0}\mathcal{A}^{i}(X)[-i],$$
and except in the very special case where $X$ is a smooth scheme, this  infinite product does not 
restrict to a finite product in general (i.e. $\mathcal{A}^{p}(X)\neq 0$ for arbitrary large $p$'s in general). 
The complex $\mathcal{A}(X)$ can be shown to carry an extra differential called the \emph{de Rham differential} and
denoted by $dR$\footnote{The existence of this differential is not a trivial fact because of the stackyness of $X$, see \cite{ptvv}.}. 
The differential $dR$ commutes with the cohomological differential, and
the complex $\mathcal{A}(X)$ will be always considered endowed with the corresponding total differential. 
When $X$ is a smooth scheme $\mathcal{A}(X)$ is simply the algebraic de Rham complex of $X$. When $X$ is 
a singular scheme $\mathcal{A}(X)$ is the derived de Rham complex of $X$, which is known to 
compute the algebraic de Rham cohomology of $X$ (see \cite{ba}). When $X$ is a derived algebraic stack the complex
$\mathcal{A}(X)$ is by definition the (algebraic and derived) de Rham complex of $X$, and it can be shown to compute the
algebraic de Rham cohomology of the underlying algebraic stack, and thus the Betti cohomology of 
its geometric realization when $k=\mathbb{C}$. 

The de Rham complex comes equipped with a standard Hodge filtration, which is a decreasing 
sequence of sub-complexes $F^{p}\mathcal{A}(X) \subset F^{p-1}\mathcal{A}(X) \subset \mathcal{A}(X)$, 
where $F^{p}\mathcal{A}(X)$ consists of the sub-complex $\prod_{i\geq p}\mathcal{A}^{i}(X)[-i] \subset 
\prod_{i\geq 0}\mathcal{A}^{i}(X)$. The complex $F^{p}\mathcal{A}(X)[p]$ is also denoted by $\mathcal{A}^{p,cl}(X)$
and is by definition the complex of closed $p$-forms on $X$. We note here that an $n$-cocycle in 
$\mathcal{A}^{p,cl}(X)$
consists of a formal series $\sum_{i\geq 0}\omega_{i}\cdot t^{i}$, where $\omega_{i}$ an element of degree $n-i$ 
in $\mathcal{A}^{p+i}(X)$, and satisfies the infinite number of equations
$$dR(\omega_{i-1}) + d(\omega_{i}) = 0 \qquad \forall \; i\geq 0,$$
where $dR$ is the de Rham differential, $d$ is the cohomological differential and $\omega_i$ is declared to be $0$ when $i<0$. With this notation, $\omega_0$ is the underlying $p$-form and 
the higher forms $\omega_i$ are the \emph{closeness structures}, reflecting that $\omega_0$ \emph{is closed
up to homotopy}. 

By definition a closed $p$-form of degree $n$ on $X$ is an element in $H^{n}(\mathcal{A}^{p,cl}(X))$. Any 
closed $p$-form $ \sum_{i\geq 0}\omega_{i}\cdot t^{i}$ of degree $n$ has
an underlying $p$-form $\omega_0$ of degree $n$, 
and thus defines a morphism $\wedge_{\OO_{X}}^{p}\mathbb{T}_{X} \longrightarrow \OO_{X}[n]$
in $L_{qcoh}(X)$. We note here that a given $p$-form of degree $n$ can come from many different 
closed $p$-forms of degree $n$, or in other words that the projection map $\mathcal{A}^{p,cl}(X) \longrightarrow 
\mathcal{A}^{p}(X)$, sending $\sum_{i\geq 0}\omega_{i}\cdot t^i$ to $\omega_{0}$, needs not be injective in cohomology. This aspect
presents  a major difference with the setting of differential forms on smooth schemes, for which 
a given $p$-forms is either closed or not closed. 
This aspect can also be understood in the setting of cyclic homology, as
differential forms on $X$ can be interpreted as elements in Hochschild homology of $X$ (suitably defined to encode
the eventual stackyness of $X$), and closed forms as elements in negative cyclic homology. 

\subsection{Shifted symplectic structures} 

\begin{df}{(\cite[Def. 1.18]{ptvv})}\label{dsymp}
An \emph{$n$-shifted symplectic structure on $X$} consists of a closed $2$-form of degree $n$ whose
underlying morphism
$$\wedge_{\OO_{X}}^{2}\mathbb{T}_{X} \longrightarrow \OO_{X}[n]$$
is non-degenerate: the adjoint map $\mathbb{T}_{X} \longrightarrow \mathbb{L}_{X}[n]$ is an equivalence of quasi-coherent 
complexes on $X$.
\end{df}

There are some basic examples of $n$-shifted symplectic structures which are the building
blocks of more evolved examples. 
A $0$-shifted symplectic structure on a smooth scheme
is simply a symplectic structure understood in the usual sense.  
For a reductive algebraic
group $G$, the $2$-shifted symplectic structures on the stack $BG$ are in one-to-one correspondence with 
non-degenerate and $G$-invariant scalar products on the Lie algebra $\frak{g}$ of $G$. Such a
structure always exists and is even unique up to a constant when $G$ is a simple reductive group. When $G=Gl_n$, 
there is a canonical choice for a 2-shifted symplectic structure on $BG$ by considering the standard 
invariant scalar product on the space of matrices given by $(A,B) \mapsto Tr(A.B)$. 

Another source of examples is provided by shifted cotangent bundles. 
For $X$ a derived algebraic stack and $n$
an arbitrary integer we define the $n$-shifted total cotangent derived stack of $X$ by
$$T^{*}X[n] := \mathbb{V}(\mathbb{L}_{X}[n]) = \Spec\, (Sym_{\OO_{X}}(\mathbb{T}_{X}[-n])),$$
as the linear derived algebraic stack over $X$ determined by the perfect complex $\mathbb{L}_{X}[n]$. 
The derived algebraic stack $T^*X[n]$ comes equipped with a standard Liouville $1$-form of degree $n$\footnote{This 
form represents the universal $1$-form of degree $n$ on $X$.}, whose
de Rham differential provides an $n$-shifted symplectic structure on $X$. This is already interesting for 
$X$ a smooth scheme as it provides instances of $n$-shifted symplectic structures for arbitrary values
of $n$. Note here that when $X$ is a smooth scheme, then  $T^*X[n]$ is either a smooth (and thus non-derived) 
algebraic $n$-stack if $n\geq 0$, or a derived scheme when $n<0$. Another interesting and useful example is when 
$X=BG$ and 
$n=1$, as $T^*X[1]$ is then identified with the quotient stack $[\frak{g}^*/G]$, for the co-adjoint action 
of $G$. The quotient stack $[\frak{g}^*/G]$ is thus equipped with a canonical $1$-shifted symplectic structure, 
which sheds new light on symplectic reduction
(we refer to \cite[\S 2.2]{ca} for more on the subject). A third important example is the
derived algebraic stack of perfect complexes $\Parf$ (see \cite[\S 2.3]{ptvv}), which 
is a generalization of the stack $BGl_n$ ($BGl_n$ sits as an open in $\Parf$).

More evolved 
examples of shifted symplectic structures can 
be constructed by means of the following existence theorem. This
result can be seen as a geometrico-algebraic counter part of the so-called AKSZ formalism. 

\begin{thm}{(\cite[Thm. 2.5]{ptvv})}\label{tsymp}
Let $X$ be either a connected compact oriented topological manifold of dimension $d$, or 
a connected smooth and proper scheme of dimension $d$ equipped with a nowhere vanishing top form $s \in 
\Omega^{d}_X$. Let 
$Y$ be a derived algebraic
stack endowed with an $n$-shifted symplectic structure. Then the derived algebraic stack 
$\Map(X,Y)$ is equipped with a canonical $(n-d)$-shifted symplectic structure. 
\end{thm}

An important special case is when $Y=BG$ for $G$ a reductive algebraic group, equipped with the $2$-shifted
symplectic structure corresponding to a non-degenerate element in $Sym^{2}(\frak{g}^{*})^G$. We find this way 
that 
the derived moduli stack of $G$-bundles on $X$, $Bun_{G}(X):=\Map(X,BG)$ carries a canonical $(2-d)$-shifted 
symplectic
structure, which is a first step towards a mathematical formulation of our principle \ref{princ}.

\begin{cor}\label{ctsymp}
With the above notations, $Bun_{G}(X)$  carries a canonical $(2-d)$-shifted symplectic structure. 
\end{cor}

When $d=2$ the above corollary recovers the well known symplectic structures on moduli spaces of
$G$-local systems on a compact Riemann surface and of $G$-bundles on K3 and abelian surface. However, even in 
this
case, the corollary is new and contains more as the $0$-shifted symplectic structure exists on the whole derived 
moduli stack, not
only on the nice part of this moduli stack which is a smooth scheme (see for instance our comments
in \S 6.1).

In dimension $3$ the corollary states that $Bun_{G}(X)$ is equipped with a natural $(-1)$-shifted symplectic
structure. The underlying $2$-form of degree $-1$ is an equivalence of perfect complexes
$\mathbb{T}_{Bun_{G}(X)} \simeq \mathbb{L}_{Bun_{G}(X)}[-1]$. When restricted to the underived part of the moduli
stack this equivalence recovers the symmetric obstruction theory used in Donaldson-Thomas theory 
(see \cite[Def. 1.1]{bu}). 
However, here again the
full data of the $(-1)$-shifted symplectic structure contains strictly more than the underlying 
symmetric obstruction theory, essentially because of the fact that a shifted symplectic structure
is not uniquely determined by its underlying $2$-form (see \cite{pato}).

Finally, when the dimension $d$ is different from $2$ and $3$ the content of the corollary seems completely new, 
thought 
in dimension $1$ it essentially states that $[G/G]$ is $1$-shifted symplectic, which can be used in order to 
provide a new understanding of quasi-hamiltonian actions (see \cite[\S 2.2]{ca}). 

The idea of the proof of theorem \ref{tsymp} is rather simple, and at least the underlying $2$-form can be 
described explicitly in terms of the formula for the tangent complexes (formula \eqref{tangent} of \S 3.4).
We define a pairing of degree $(n-d)$ on this complex by the composition
of the natural maps and the pairing of degree $n$ on $\mathbb{T}_{Y}$
$$\xymatrix{
\wedge^{2}p_{*}(ev^{*}(\mathbb{T}_{Y}))  \ar[r] & 
p_{*}(ev^{*}(\wedge^{2}\mathbb{T}_{Y})) \ar[r] & 
p_{*}(\OO)[n] \ar[r] & \OO[n-d],}$$
where the last map comes from the fundamental class in $H^{d}(X,\OO_{X})\simeq k$. This defines a non-degenerate 
$2$-form on 
$\Map(X,Y)$, and the main content of the theorem \ref{tsymp} is that this $2$-form comes from a canonically 
defined
closed $2$-form of degree $(n-d)$. 
 
For variants and generalizations of theorem \ref{tsymp} we refer to 
\cite{ptvv,ca,dagems} in which the reader will find non-commutative generalizations as well 
as versions with boundary conditions, but also several other possible admissible sources. 

\subsection{Derived critical loci}
To finish the part on $n$-shifted symplectic structures let us mention critical loci and their 
possible generalizations. 
We have already seen that for a given derived algebraic stack $X$ the shifted cotangent $T^*X[n]$ carries a 
canonical
$n$-shifted symplectic structure. Moreover, the zero section $X \longrightarrow T^*X[n]$ has a natural Lagrangian 
structure (see \cite[Def. 2.8]{ptvv}). 
More generally, if $f \in H^{n}(X,\OO_{X})$ is a function of degree $n$ on $X$, its de Rham differential $dR(f)$ 
defines a morphism of derived algebraic stacks $dR(f) : X \longrightarrow T^{*}X[n]$
which is also equipped with a natural Lagrangian structure. Therefore, the intersection of the zero section with the
section $dR(f)$ defines a natural $(n-1)$-shifted derived algebraic stack 
(see \cite[Thm. 2.9]{ptvv}) denoted by $\mathbb{R}Crit(f)$ and called
the derived critical locus of $f$. When $f=0$ the derived critical locus $\mathbb{R}Crit(f)$ is simply 
$T^{*}X[n-1]$ together
with its natural $(n-1)$-shifted symplectic structure. When $X$ is a smooth scheme and $f$ is a function of 
degree $0$ (i.e. simply a function $X \longrightarrow \mathbb{A}^{1}$), then the symplectic geometry of 
$\mathbb{R}Crit(f)$ is closely related to the singularity theory of the function $f$. 
From a general point of view derived critical loci provide a nice source of examples of $n$-shifted symplectic derived algebraic
stacks, which contain already examples of geometric interests. It is shown in \cite{brbu,bogr} 
that every $(-1)$-shifted symplectic
derived scheme is locally the derived critical locus of a function defined on a smooth scheme. 

Derived critical loci are important because they are easy to describe and their quantizations
can be understood explicitly. Moreover, derived critical loci and their
generalizations can be used to provide local models for $n$-shifted symplectic
structures by means of a formal Darboux lemma we will not reproduce here (see for instance
\cite{bbbj,brbu,bogr}).

\subsection{Shifted polyvector fields and poisson structures} The notion of shifted Poisson structure
 is the dual notion of that of shifted 
symplectic structure we have discussed so far.  The general theory of shifted Poisson structures has not been 
fully settled down yet and we will here present the basic definitions as well  as its, still hypothetical, 
relations with shifted symplectic structures. They are however a key notion in the existence of quantization
that will be presented in the next section. 

\subsubsection{Shifted polyvectors on derived algebraic stacks} A derived algebraic stack $X$ (as usual assumed 
locally of finite presentation over the ground field $k$) 
has a tangent complex $\mathbb{T}_{X}$, which is the $\OO_{X}$-linear dual to the cotangent complex. 
The complex of $n$-shifted polyvector fields on $X$ is defined by
$$\displaystyle{\Pol(X,n) := \bigoplus_{i} \Gamma(X,Sym^{i}_{\OO_{X}}(\mathbb{T}_{X}[-1-n]))}.$$
The complex $\Pol(X,n)$ has a natural structure of a graded commutative dg-algebra, for which the piece of 
weight $i$ is $\Gamma(X,Sym^{i}_{\OO_{X}}(\mathbb{T}_{X}[-1-n]))$ and the multiplication is induced by the 
canonical multiplication 
on the symmetric algebra. We note here that depending of the parity of $n$ we either
have $Sym^{i}_{\OO_{X}}(\mathbb{T}_{X}[-1-n]) \simeq (\wedge^{i}\mathbb{T}_{X})[-i-ni]$ (if $n$ is even), or
$Sym^{i}_{\OO_{X}}(\mathbb{T}_{X}[-1-n]) \simeq (Sym^{i}\mathbb{T}_{X})[-i-ni]$ (if $n$ is odd).
When $X$ is a smooth scheme and $n=0$, $\Pol(X,0)=\oplus_{i}\Gamma(X,\wedge^{i}\mathbb{T}_{X})[-i]$ is the 
standard
complex of polyvector fields of $X$. When $n=1$, and still $X$ a smooth scheme, 
$\Pol(X,1)$ coincides with $\Gamma(T^*X,\OO_{T^*X})$, the cohomology of the total cotangent space of
$X$ with coefficients in $\OO$. In general, $\Pol(X,n)$ can be interpreted as the graded cdga of cohomology of
the shifted cotangent derived stack $T^*X[n+1]$ with coefficients in $\OO$ (that is "functions" on $T^*X[n+1]$). 

When $X$ is a smooth scheme, $\mathbb{T}_{X}$ is a sheaf (say on the small \'etale site of $X$) 
of $k$-linear Lie algebras with the bracket of vector fields. 
This extends easily to the case where $X$ is a derived Deligne-Mumford stack, 
$\mathbb{T}_{X}$ can be made into a sheaf of $k$-linear dg-Lie algebras for the bracket of dg-derivations. 
Therefore, 
polyvector fields $\Pol(X,n)$ can also be endowed with a $k$-linear dg-Lie bracket of cohomological degree 
$-1-n$, making it into
a graded Poisson dg-algebra where the bracket has cohomological degree $(-1-n)$ and weight $(-1)$. In particular 
$\Pol(X,n)[n+1]$ always comes equipped with a structure of a graded dg-Lie algebra over $k$. It is expected that 
this fact 
remains valid for a general derived algebraic stack $X$, but there is no precise construction at the moment. One 
complication
when considering general algebraic stacks comes from the fact that vector fields can not be pulled-back along 
smooth morphisms
(as opposed to \'etale maps), making the construction of the Lie 
bracket on $\Pol(X,n)$ much more complicated
than for the case of a scheme. For a derived algebraic stack of the form $[\Spec\, A/G]$, for $G$ linear, 
there are however two possible constructions. A first very indirect construction uses 
natural operations  
on the derived moduli stacks of branes (see \cite{to3}). A more direct construction can be 
done as follows. We can take $A$ to be a cofibrant and fibrant cgda inside the category of representations 
$Rep(G)$. 
We let $\mathbb{T}_{A}$ be the $A$-dg-module of dg-derivations from $A$ to itself. The action of $G$ on $A$ 
induces
a morphism of dg-Lie algebras $\frak{g}\otimes_{k} A \longrightarrow \mathbb{T}_{A}$ representing the 
infinitesimal action
of $G$ on $A$.
We consider 
the co-c\^one $\mathbb{T}$ of the morphism $\frak{g}\otimes_{k} A \longrightarrow \mathbb{T}_{A}$. The 
complex
$\mathbb{T}$ is obviously a $k$-linear Lie algebra for the bracket induced from the brackets on 
$\mathbb{T}_{A}$
 and on $\frak{g}$, but this lie structure is \emph{not} compatible with the cohomological differential and 
 thus
is not a dg-Lie algebra. However, its fixed points by $G$ (assume $G$ reductive for simplicity) 
is a dg-Lie algebra over $k$, which is a model for $\Gamma(X,\mathbb{T}_{X})$ where $X=[\Spec\, A/G]$. This
construction can be also applied to the $G$-invariant of the various symmetric powers of shifts of 
$\mathbb{T}$ in
order to get the desired dg-Lie structure on $\Pol(X,n)[n+1]$ in this special case. \vspace{8mm}

\subsubsection{Shifted Poisson structures}

Let $X$ be a derived algebraic stack and fix an integer $n\in \mathbb{Z}$.  We can define $n$-shifted 
Poisson structures
as follows. We let $\Pol(X,n)[n+1]$ be the shifted polyvector fields on $X$, endowed with the structure of a 
graded 
dg-Lie algebra just mentioned. We let $k(2)[-1]$ be the graded dg-Lie algebra which is $k$ in cohomological 
degree $1$, 
with zero bracket and $k$ is pure of weight $2$. An $n$-shifted Poisson structure on $X$ is then 
defined to be a morphism of graded dg-Lie algebras
$$p : k(2)[-1] \longrightarrow \Pol(X,n)[n+1].$$
Here, a morphism of graded dg-Lie algebras truly means a morphism inside the $\s$-category of
graded dg-Lie algebras, or a morphism in an appropriate homotopy category. Using the
dictionary between dg-Lie algebras and formal moduli problems (see \cite{lu2}), such a morphism $p$ is 
determined by 
a Mauer-Cartan element in $\Pol(X,n)[n+1] \otimes tk[[t]]$, which is of weight $2$ with respect to the 
grading
on $\Pol(X,n)$. Such an element can be described explicitly as a formal power series 
$\sum_{i\geq 1}p_{i}\cdot t^i$, where $p_{i}$ is an element of cohomological degree $n+2$ in 
$\Gamma(X,Sym^{i+1}(\mathbb{T}_{X}[-1-n]))$, and satisfies the equations
$$\displaystyle{d(p_i) + \frac{1}{2}\cdot \sum_{a+b=i}[p_a,p_b] = 0 \qquad \forall i\geq 1}.$$

As we already  mentioned, shifted Poisson structures can be developed along the same lines as
shifted symplectic structures (e.g. there is a notion of co-isotropic structures on a map with an $n$-shifted 
Poisson target, 
and a Poisson version of the existence theorem \ref{tsymp}), 
but at the moment this work has not been carried out in full details. It is believed that 
for a given $X$ and $n\in \mathbb{Z}$, there is a one-to-one correspondence between $n$-shifted symplectic 
structures
on $X$ and $n$-shifted Poisson structures on $X$ which are non-degenerate in an obvious sense. However, this
correspondence has not been established yet, except in some special cases, and remains at the moment an open 
question
for further research (see \S 6.2). \\

\section{Deformation quantization of $n$-shifted Poisson structures}

In this section we finally discuss the existence of quantization of $n$-shifted Poisson structures, a far 
reaching 
generalization of the existence of deformation quantization of Poisson manifolds due to Kontsevich. For this
we first briefly discuss the output of the quantization, namely the notion of deformation of categories and 
iterated monoidal
categories, which already contains some non-trivial aspects. We then present the formality conjecture, which 
is
now a theorem except in some very particular cases,  and whose main corollary is the fact that every 
$n$-shifted 
Poisson structure
defines a canonical formal deformation of the $E_n$-monoidal category of quasi-coherent complexes. We also 
discuss 
the case $n<0$ by presenting the \emph{red shift trick} consisting of working with a formal parameter $\hbar$ 
living in
some non-trivial cohomological degree.

\subsection{The deformation theory of monoidal dg-categories}

As we have seen in \S 2, a derived algebraic stack $X$ has a dg-category of quasi-coherent complexes $L(X)$. 
It is a $k$-linear dg-category which 
admits arbitrary colimits. We will assume in this section that $L(X)$
is a compactly generated dg-category, or equivalently that it can be realized as the category of 
dg-modules over a small 
dg-category. More generally we will assume that $X$ is a \emph{perfect} derived algebraic stack,
in the sense that perfect complexes on $X$ are compact generators of $L(X)$. This is known to be the case 
under the assumption that $X$ can be written as a quotient $[\Spec\, A/G]$ for a linear algebraic $G$ acting
on a cdga $A$. \\

\subsubsection{Deformations of dg-categories}
We let $T_0:=L(X)$ and we would like to study the deformation theory of $T_0$. For this, we define a first naive 
deformation functor
$Def^{naive}(T_0) : \art \longrightarrow \Top$, from the $\s$-category of augmented local artinian cdga to the $\s$-
category 
of spaces as follows. To $A\in \art$ we assign the $\s$-category $\mathbb{D}g^{c}(A)$, of cocomplete and
compactly generated $A$-linear dg-categories  and $A$-linear colimit preserving dg-functors
(see \cite[\S 3.1]{to4}). For a morphism of dg-artinian rings $A \rightarrow B$, we have a base change $\s$-functor
$-\widehat{\otimes}_{A}B : \mathbb{D}g^{c}(A) \longrightarrow \mathbb{D}g^{c}(B)$.  We then set 
$$Def^{naive}(T_0)(A) := \mathbb{D}g^{c}(A) \times_{\mathbb{D}g^{c}(k)}\{T_{0}\}.$$
Here $\mathbb{D}g^{c}(A) \times_{\mathbb{D}g^{c}(k)}\{T_{0}\}$ is the fiber taken at the point $T_0$ of the 
$\s$-functor $-\widehat{\otimes}_{A}k$ induced by the augmentation $A \rightarrow k$. As is, $Def^{naive}(T_0)(A)$ 
is an $\s$-category, 
from which we extract a space by taking the geometric realization of its sub-$\s$-category of equivalences (i.e. 
taking the
nerve of the maximal sub-$\s$-groupoid). Intuitively, $Def^{naive}(T_0)(A)$ is the classifying space of pairs 
$(T,u)$, with
$T$ a compactly generated $A$-linear dg-category and $u$ a $k$-linear equivalence $u : T\widehat{\otimes}_{A}k 
\simeq T_0$. 

As already observed in \cite{kelo} the $\s$-functor $Def^{naive}(T_0)$ is not a formal moduli problem, 
it does not satisfies the Schlessinger conditions of \cite{lu2}, and thus
can not be equivalent to the functor of Mauer-Cartan elements in a dg-Lie algebra. This bad behavior of 
the $\s$-functor $Def^{naive}(T_0)$ has been a longstanding major obstacle preventing the
 understanding of the deformation theory 
of dg-categories. There have been several tentative modifications of $Def^{naive}(T_0)$ attempting to overcome this 
problem, for 
instance by allowing curved dg-categories as possible deformations, however none of these were 
successful. We propose here a new 
solution to 
this problem
which provides the only complete understanding of deformations of dg-categories that we are aware of.
For this, we introduce $Def(T_0) : \art \longrightarrow \Top$, which is the universal $\s$-functor constructed out of 
$Def^{naive}(T_0)$ and satisfying the Schlessinger conditions of \cite{lu2} (in other words it is the best possible
approximation of $Def^{naive}(T_0)$ by an $\s$-functor associated to a dg-Lie algebra). By construction there is a 
natural transformation $l : Def^{naive}(T_0) \longrightarrow Def(T_0)$, as well as a dg-Lie algebra $L$ such that 
$Def(T_0)$ is given by $A \mapsto \underline{MC}_{*}(L \otimes m_A)$ (where as usual $m_A \subset A$ is the augmentation
dg-ideal in $A$, and $\underline{MC}_{*}$ denotes the space of Mauer-Cartan elements). Moreover the natural 
transformation $l$ is universal for these properties, and in particular the dg-Lie $L$ is uniquely determined 
and only depends on $Def^{naive}(T_0)$. 

The following theorem is folklore and known to experts. It appears for instance 
in a disguised form in \cite{pre}.

\begin{thm}\label{tdef}
Let $T_0$ be a compactly generated dg-category. 
\begin{enumerate}
\item The dg-Lie algebra associated to the formal moduli problem $Def(T_0)$ is $HH(T_0)[1]$, the
Hochschild cochains on $T_0$ endowed with its usual Gerstenhaber bracket (see e.g. \cite[\S 5.4]{ke}).

\item The space $Def(T_0)(k[[t]])$ is naturally equivalent to 
the classifying space of $k[\beta]$-linear structures on $T_0$, where $k[\beta]$
is the polynomial dg-algebra over $k$ with one generator $\beta$ in degree $2$. 

\end{enumerate}
\end{thm}

The above theorem subsumes the two main properties of the formal moduli problem $Def(T_0)$,
but much more can be said. The formla for the $k[[t]]$-points of $Def(T_0)$ can be generalized to 
any (pro-)artinian augmented dg-algebra $A$, by using $B_{A}$-linear structures on $T_0$, where
now $B_A$ is the $E_2$-Koszul dual of $A$ (see \cite{lu2}, and the $E_2$-Koszul dual of $k[[t]]$ is of course
$k[\beta]$). By construction we have a map of spaces
$Def^{naive}(T_0)(A) \longrightarrow Def(T_0)(A)$. This map is not an equivalence but can be shown to 
have $0$-truncated fibers (so it induces isomorphisms on $\pi_i$ for $i>1$ and is injective on $\pi_1$). 
It is interesting to note here that not only $Def(T_0)$ contains more objects than $Def^{naive}(T_0)$
but also contains more morphisms. There are natural conditions one can impose on $T_0$ in order to make
$Def^{naive}(T_0)$ closer to $Def(T_0)$. It is for instance believed that they coincide when 
$T_0$ is a smooth and proper dg-category, as well as for dg-categories of complexes in Grothendieck abelian categories.
In our situation, $T_0=L(X)$, with $X$ a derived algebraic stack which is not smooth in general, it is not reasonable
to expect any nice assumptions on $T_0$, and the above theorem is probably the best available result in order
to understand formal deformations of $L(X)$.  \\

\subsubsection{Deformations of monoidal dg-categories}

Theorem \ref{tdef} also possesses monoidal and iterated monoidal versions as follows. First of all 
the $\s$-category $\mathbb{D}g^{c}(A)$ of compactly generated $A$-linear dg-categories is equipped with a 
tensor product $\widehat{\otimes}_{A}$, making it into 
a symmetric monoidal $\s$-category. It is therefore possible to use the notion
of an $E_n$-monoid in $\mathbb{D}g^{c}(A)$ of \cite{ha}, in order to define
$E_n$-monoidal $A$-linear dg-categories (also called $n$-fold monoidal $A$-linear dg-categories). In a nutshell,
an $E_n$-monoidal $A$-linear dg-category consists of a compactly generated $A$-linear dg-category $T$
together with morphisms 
$$\mu_k : E_n(k) \otimes T^{\hat{\otimes_A}\, k} \longrightarrow T,$$
where the tensor by the space $E_n(k)$ and the tensor 
products are taken in the symmetric monoidal $\s$-category $\mathbb{D}g^{c}(A)$, and
together with compatibilty conditions/structures.
For 
our derived algebraic stack $X$, the dg-category $L(X)$ is equipped with a symmetric monoidal 
structure and thus is naturally an $E_n$-monoidal dg-category for all $n\geq 0$, where by convention an 
$E_0$-monoidal dg-category simply is a dg-category. \\

For a cdga $A$, we set $\mathbb{D}g^{c}_{E_{n}}(A)$ for the $\s$-category of compactly generated 
$E_n$-monoidal $A$-linear dg-categories. Here compactly generated also means that the compact objects are
stable by the monoidal structure, so objects in $\mathbb{D}g^{c}_{E_{n}}(A)$ can also be described 
as dg-categories of dg-modules over small $A$-linear $E_n$-monoidal dg-categories. Morphisms in 
$\mathbb{D}g^{c}_{E_n}(A)$ must be defined with some care as they involve higher dimensional 
versions of Morita morphisms between algebras. For two $E_n$-monoidal dg-categories $T$ and $T'$ in 
$\mathbb{D}g^{c}_{E_n}(A)$, the dg-category of $A$-linear colimit preserving dg-functors
can be written as $T^{\vee} \widehat{\otimes}_{A}T'$, where $T^{\vee}$ is the dual of $T$ (i.e. we take
the opposite of the sub-dg-category of compact generators). The $A$-linear dg-category $T^{\vee} \widehat{\otimes}_{A}T'$
is a new object in $\mathbb{D}g^{c}_{E_n}(A)$, and in particular it makes sense to consider $E_n$-algebras
inside the dg-category $T^{\vee} \widehat{\otimes}_{A}T'$. From the point of view of dg-functors these correspond to 
$E_n$-lax monoidal $A$-linear colimit preserving dg-functors $T \rightarrow T'$. For two 
$E_n$-algebras $M$ and $N$ inside $T^{\vee} \widehat{\otimes}_{A}T'$, we can form a new 
$E_n$-algebra $M^{op}\otimes N$. The $A$-linear dg-category of $M^{op}\otimes N$-modules inside
$T^{\vee} \widehat{\otimes}_{A}T'$ is then $E_{n-1}$-monoidal, so the process can be iterated. We can 
consider two $E_{n-1}$-algebras inside $M^{op}\otimes N$-modules, say $M'$ and $N'$, as well
as their tensor product $M'^{op}\otimes N'$ and the $A$-linear dg-category of 
$M'^{op}\otimes N'$-modules, which is itself $E_{n-2}$-monoidal \dots and so on and so forth. We are describing 
here $\mathbb{D}g^{c}_{E_{n}}(A)$ as an $(\s,n+1)$-category (see \cite{si1}), whose objects are $E_n$-monoidal 
compactly generated $A$-linear dg-categories, whose $1$-morphism from $T$ to $T'$ are $E_{n-1}$-algebras inside
$T^{\vee} \widehat{\otimes}_{A}T'$, whose $2$-morphisms between $M'$ and $N'$ are $E_{n-2}$-algebras inside
$M'^{op}\otimes N'$-modules, etc \dots. The $(\s,n+1)$-category $\mathbb{D}g^{c}_{E_n}(A)$ produces a space 
by considering the geometric realization of its maximal sub-$\s$-groupoid (i.e. realizing the sub-$\s$-category 
of equivalences). \\

For $T_0=L(X)$, assuming that $L(X)$ is compactly generated and that its compact objects are the perfect complexes, 
we define a naive deformation functor $Def^{naive}_{E_{n}}(T_0)$, of $T_0$ considered as an $E_n$-dg-category, by 
sending an augmented dg-artinian ring $A \in \art$ to 
the fiber at $T_0$ of the restriction map
$$-\widehat{\otimes}_{A} k : \mathbb{D}g^{c}_{E_n}(A) \longrightarrow \mathbb{D}g^{c}_{E_n}(k).$$
The space $Def^{naive}_{E_{n}}(T_0)$ is the space of pairs $(T,u)$, where $T$ is an 
$E_n$-monoidal compactly generated
$A$-linear dg-category, and $u : T\widehat{\otimes}_{A}k \simeq T_0$ an equivalence in 
$\mathbb{D}g^{c}_{E_{n}}(k)$. Similar to the case $n=0$ we already discussed, the $\s$-functor 
$Def^{naive}_{E_{n}}(T_0)$ does not satisfy the Schlessinger's conditions, and the bigger $n$ is,
 the more this fails.
We denote by $Def_{E_n}(T_0)$ the formal moduli problem generated by $Def^{naive}_{E_{n}}(T_0)$. The following 
theorem is the generalization of \ref{tdef} to the iterated monoidal setting.

\begin{thm}\label{tdef2}
Let $T_0$ be a compactly generated $E_n$-monoidal dg-category. 
\begin{enumerate}
\item The dg-Lie algebra associated to the formal moduli problem $Def_{E_n}(T_0)$ is $HH^{E_{n+1}}(T_0)[n+1]$, the
$E_{n+1}$-Hochschild cochains on $T_0$ of \cite{fr}.

\item The space $Def_{E_n}(T_0)(k[[t]])$ is naturally equivalent to 
the classifying space of $k[\beta_{n}]$-linear structures on $T_0$, where $k[\beta_{n}]$
is the commutative polynomial dg-algebra over $k$ with one generator $\beta_n$ in degree $2+n$. 

\end{enumerate}
\end{thm}

Theorems \ref{tdef} and \ref{tdef2} provides a way to understand the relations between (higher) Hochschild 
cohomology and
formal deformations of dg-categories and iterated monoidal dg-categories. They state in particular that the correct 
manner
to define a formal deformation of a given dg-category $T_0$, parametrized by $k[[t]]$, is by 
considering $k[\beta]$-linear structures on $T_0$, and similarly for the iterated monoidal setting 
with $k[\beta_n]$-linear structures. In the sequel, we will freely use the expression "formal deformation of 
the dg-category $L(X)$ considered as an $E_n$-monoidal dg-category", by which we mean 
an element in $Def_{E_{n}}(L(X))(k[[t]])$, and thus a $k[\beta_{n}]$-linear structure on $T_0$. We however continue 
to think of these deformations as actual deformations of $L(X)$ over $k[[\hbar]]$ for a formal parameter $\hbar$, 
even thought 
they are not quite as naive objects. 

\subsection{The higher formality conjecture}
We have just seen that formal deformations of a given $E_n$-monoidal compactly generated dg-category $T_0$ is controlled
by its higher Hochschild cochain complex $HH^{E_{n+1}}(T_0)[n+1]$, endowed with its natural structure of a dg-Lie algebra.
We now turn to the specific case where $T_0=L(X)$, the quasi-coherent dg-category of a derived algebraic stack
$X$. We continue to assume that $X$ is nice enough (e.g. of the form $[\Spec\, A/G]$) so that $L(X)$ is
compactly generated by the perfect complexes). The higher Hochschild cohomology of $T_0$ can then be described in geometric
terms as follows. We let $n\geq 0$, and let $S^{n}=\partial B^{n+1}$ 
be the topological  $n$-sphere considered as a constant derived stack.
We consider the derived mapping stack $\mathcal{L}^{(n)}(X):=\Map(S^n,X)$, also called the 
\emph{$n$-dimensional derived loop stack of $X$}.
There is a constant map morphism $j : X \longrightarrow \mathcal{L}^{(n)}(X)$, and thus a quasi-coherent 
complex $j_{*}(\OO_{X}) \in L(\mathcal{L}^{(n)}(X))$. The $E_{n+1}$-Hochschild cohomology of the dg-category $L(X)$ can
be identified with 
$$HH^{E_{n+1}}(X)\simeq End_{L(\mathcal{L}^{(n)}(X))}(j_{*}(\OO_X)).$$
Note that when $n=0$ and $X$ is a scheme this recovers the description of the Hochschild complex of $X$ as the self extension
of the diagonal. Because of the stackyness of $X$ this definition can be modified by replacing $\mathcal{L}^{(n)}(X)$
by its formal completion $\widehat{\mathcal{L}^{(n)}(X)}$ along the map $X \longrightarrow \mathcal{L}^{(n)}(X)$, which 
is called the formal $n$-dimensional derived loop space (when $X$ is a derived scheme the formal and non-formal 
versions of the derived loop stacks coincide). We then have the formal version of Hochschild complex
$$\widehat{HH}^{E_{n+1}}(X)\simeq End_{L(\widehat{\mathcal{L}^{(n)}(X)})}(j_{*}(\OO_X)).$$
Note that we have a natural morphism $\widehat{HH}^{E_{n+1}}(X) \longrightarrow HH^{E_{n+1}}(X)$. 

The complexe $\widehat{HH}^{E_{n+1}}(X)$ has a structure of an $E_{n+2}$-algebra, as predicted by the so-called Deligne's conjecture
which is now a theorem (see \cite{fr,ha}). In particular $\widehat{HH}^{E_{n+1}}(X)[n+1]$ is a dg-Lie algebra. 
The formality conjecture asserts
that the dg-Lie algebra $\widehat{HH}^{E_{n+1}}(X)[n+1]$ can be described in simple terms involving shifted polyvector fields.

\begin{conj}{(Higher formality)}\label{cform}
For a nice enough derived algebraic stack $X$, and $n\geq 0$, the dg-Lie algebra $\widehat{HH}^{E_{n+1}}(X)[n+1]$
is quasi-isomorphic to $\Pol(X,n)[n+1]$. The quasi-isomorphism is canonical up to a universal choice of a Drinfeld
associator.
\end{conj}

Note that when $X$ is a smooth scheme and $n=0$ the conjecture \ref{cform} 
is the so-called \emph{Kontsevich's formality theorem}. 
The conjecture has been proven in already many cases.

\begin{thm}{(\cite[Cor. 5.4]{to3})}\label{tform}
The above higher formality conjecture is true for all $n>0$ and for all derived algebraic stacks $X$ of the form
$[\Spec\, A/G]$ for $G$ a linear algebraic group acting on the cdga $A$. When $X$ is a derived Deligne-Mumford stack it is
also true for $n=0$.
\end{thm}

The theorem above provides many cases in which conjecture \ref{cform} is satified. We believe it is also true in the remaining
case when $n=0$ and for non Deligne-Mumford stacks. We also believe that the restriction 
for $X$ being of the form
$[\Spec\, A/G]$ in the theorem \ref{tform} is not necessary, and that the theorem should be 
true for a large class of derived higher algebraic stacks as well. 

\subsection{Existence of deformation quantization}
We finally arrive at the existence of quantization of derived algebraic stacks $X$ endowed with $n$-shifted 
Poisson structures, 
and its consequence: the mathematical incarnation of our principle \ref{princ}. Let $X$ be a derived algebraic stack, and $n \geq 0$ to start with (the case of negative values will be treated below). 
We assume that $X$ is nice enough and that the conjecture \ref{tform} is satisfied (e.g. under the
hypothesis of theorem \ref{tform}). Let $p$ be an $n$-shifted Poisson structure
on $X$. By definition it provides a morphism of  dg-Lie algebras
$p : k[-1] \longrightarrow Pol(X,n+1)[n+1].$
Using the conjecture \ref{tform} we find a morphism of dg-Lie algebras
$p : k[-1] \longrightarrow \widehat{HH}^{E_{n+1}}(X)$, which composed with the natural morphism 
$\widehat{HH}^{E_{n+1}}(X) \longrightarrow
HH^{E_{n+1}}(X)$ provides a morphism of dg-Lie algebras
$$p : k[-1] \longrightarrow HH^{E_{n+1}}(X).$$
The derived deformation theory (see \cite{lu2}) and theorem \ref{tdef2} tell 
us that the morphism $p$ provides a formal
deformation of $L(X)$ as an $E_n$-monoidal dg-category, denoted by $L(X,p)$. This is the deformation quantization 
of the pair $(X,p)$. 

Assume now that $n<0$ and that $X$ is equipped with an $n$-shifted Poisson structure $p$ such that the
conjecture \ref{tform} is satisfied for $X$ and $-n$. The $n$-shifted Poisson structure $p$
is a morphism of graded dg-Lie algebras
$k(2)[-1] \longrightarrow \Pol(X,n)[n+1]$
where $k(2)[-1]$ is the abelian dg-Lie algebra which is $k$ in cohomological degree $1$ and pure weight $2$. 
The category of $\mathbb{Z}$-graded complexes has a tensor auto-equivalence, sending a complex
$E$ pure of weight $i$ to $E[-2i]$ again pure of weight $i$. This auto-equivalence induces an auto-equivalence of the
$\s$-category of graded dg-Lie algebras, and sends $\Pol(X,n)[n+1]$ to $\Pol(X,n+2)[n+3]$
and $k(2)[-1]$ to $k(2)[-3]$. Iterated $n$ times, the morphism $p$ goes to a new morphism of dg-Lie algebras
$p' : k(2)[-2n-1] \longrightarrow \Pol(X,-n)[-n+1]$, which by conjecture \ref{cform} induces a morphism 
of dg-Lie algebras
$$p' : k[-2n-1] \longrightarrow HH^{E_{-n+1}}(X)[-n+1].$$
The abelian dg-Lie algebra $k[-2n-1]$ corresponds to the formal derived scheme $\Spec\, k[[\hbar_{2n}]]$, 
where now $\hbar_{2n}$ has cohomological degree $2n$. By the general DDT and theorem \ref{tdef2} we do find 
a formal deformation of $L(X)$, considered as an $E_{-n}$-monoidal dg-category, over $k[[\hbar_{2n}]]$. This
deformation will be denoted by $L(X,p)$. 
This trick to deal with cases where $n<0$ is called the \emph{red shift trick}. It is not new, 
and already appears in the conjecture \cite[Page 14]{ka} where $\mathbb{Z}/2$-graded derived categories are 
considered instead of $\mathbb{Z}$-graded derived categories, and canceling out the red shift.

\begin{df}\label{dquant}
The formal deformation $L(X,p)$ constructed above is the \emph{deformation quantization of $(X,p)$}. It is a
formal deformation of $L(X)$ considered as an $E_n$-monoidal dg-category if $n\geq 0$, and 
a formal deformation of $L(X)$  considered as an $E_{-n}$-dg-category over $k[[\hbar_{2n}]]$ if $n<0$. 
\end{df}

Definition \ref{dquant} applies in particular to the case $X=Bun_G(Y)$, making 
our principle \ref{princ} into a mathematical statement.

\section{Examples and open questions}

We present here some examples as well as some further questions. 
\subsection{Three examples}

We start by 
coming back to the three situations we mentioned in  \S 2. 

\textbf{Quantum groups.} We let $X=BG$, for $G$ reductive. We have seen that $X$ has a $2$-shifted symplectic 
structure given by the choice of non-degenerate $G$-invariant scalar product on $\frak{g}$. The dg-category
$L(X)$  here is the dg-category of complexes of representations of $G$. Our quantization is then 
a formal deformation of $L(X)$ as an $E_2$-monoidal dg-category, and is simply realized by taking the dg-category 
of complexes of representations of the quantum group. 

\textbf{Skein dg-algebras.} We now let $X=Bun_G(\Sigma)$ be the derived moduli stack of $G$-bundles on 
a compact oriented surface $\Sigma$. We know that $X$ carries a natural $0$-shifted symplectic structure (depending
on a choice of a non-degenerate $G$-invariant scalar product on $\frak{g}$),
whose quantization $L(X,p)$ in our sense is a deformation of the dg-category $L(X)$\footnote{Note however that here 
$n=0$ and
the formality conjecture \ref{cform} is not established yet, so this situation is still conjectural at the moment.}.
The dg-category $L(X,p)$ is an interesting refinement of the skein algebra of $\Sigma$ which, as far as the author 
is aware, 
has not been considered before. The structure sheaf $\OO_{X} \in L(X)$ deforms to a uniquely defined
object $\widetilde{\OO_{X}} \in L(X,p)$, whose endormophisms 
form a dg-algebra $B_{\hbar}=End(\widetilde{\OO_{X}})$ 
over $k[[\hbar]]$, which is
a deformation of $\OO_{X}(X)$ the dg-algebra of functions on $X$. The skein algebra is recovered as 
$H^{0}(B_{\hbar})$, but
$B_{\hbar}$ is not cohomologically concentrated in degree $0$ in general and contains strictly more than 
$K_{\hbar}(\Sigma)$.
The higher cohomology groups of $B_{\hbar}$ are directly related to the non-trivial
derived structure of $X$, which is 
concentrated around the singular points corresponding to $G$-bundles with many automorphisms. 
Outside these bad points the dg-category $L(X,p)$ is essentially given by complexes of $K_{\hbar}(\Sigma)$-modules.
Formally around a given singular point $\rho \in X$, the dg-category $L(X,p)$ has a rather simple description
as follows. The formal completion of $X$ at $\rho$ is controlled by the formal dg-Lie algebra 
$L_{\rho}:=H^{*}(\Sigma,\frak{g}_{\rho})$, where $\frak{g}_{\rho}$ is the local system of Lie algebras associated 
to the $G$-bundle $\rho$. The dg-Lie algebra $L_{\rho}$ is endowed with a non-degenerate pairing of degree $2$ 
induced by 
the choice of a $G$-invariant scalar product on $\frak{g}$ which defines a non-degenerate pairing
$p : L_{\rho}^{\vee}[-1] \wedge L_{\rho}^{\vee}[-1] \longrightarrow k$. The pairing $p$ defines itself 
a Poisson structure on the completed Chevalley complex 
$\widehat{\OO_{X,\rho}}\simeq \widehat{Sym}_{k}(L_{\rho}^{\vee}[-1])$, which is 
the cdga of formal functions on $X$ around $\rho$. The quantization of this Poisson cdga, which 
can be described in simple terms as the Weyl dg-algebra associated to $L_{\rho}$ with the pairing $p$, is the quantization 
of $X$ around $\rho$ and can be used to describe the full sub-dg-category of $L(X,p)$ 
generated by objects supported at $\rho$. 

\textbf{Donaldson-Thomas theory.} We now turn to the case where $X=Bun_G(Y)$ for a Calabi-Yau 3-fold $Y$, 
which 
is endowed with a $(-1)$-shifted symplectic form. Our quantization
$L(X,p)$ here is a formal deformation of $L(X)$ as a monoidal dg-category with a formal parameter $
\hbar_{-2}$ of degree
$-2$. To simplify a bit we can consider this as a formal deformation of $L^{\mathbb{Z}/2}(X)$, 
the $2$-periodic dg-category of quasi-coherent complexes on $X$, considered as a monoidal dg-category and
with a formal parameter $\hbar$ sitting now in degree $0$. Locally, $X$ is essentially given as the critical 
locus
of a function $f$, whose category of matrix factorizations $MF(f)$ provides a natural $L(X,p)$-module (i.e. 
$MF(f)$ is enriched over the monoidal dg-category $L(X,p)$). In a precise sense, $MF(f)$ can be viewed
as an object $\mathcal{M}$ in the quantization of $X$\footnote{This is so when monoidal dg-categories are 
considered
through their $(\s,2)$-category of modules.}.  The object $\mathcal{M}$ only exists locally, but when $X$ 
is endowed with orientation data we can expect more and maybe an existence globally on $X$ (for instance,
the class of $\mathcal{M}$ in a suitable Grothendieck group has been constructed in \cite{ks}). This 
suggests a possible relation with the perverse sheaf $\mathcal{E}$ we mentioned in \S 2, as $\mathcal{E}$
should be somehow the Betti realization of the sheaf of dg-categories $\mathcal{M}$. Our quantization 
should thus refine and reinterpret some already known constructions in Donaldson-Thomas theory. 

\vspace{1cm}

\subsection{Further questions} We finish by a sample of further possible research directions. \\

\textbf{Symplectic to Poisson and formality for $n$=0.} As already mentioned in the text the precise way to obtain 
an $n$-shifted Poisson structure out of an $n$-shifted symplectic structure is not clear at the moment, except in 
some special case (e.g. for derived scheme for which a version of the Darboux lemma holds and can be used, see \cite{brbu,bogr}).
Also recall that our conjecture \ref{cform} remains open for non Deligne-Mumford derived algebraic stacks.

\textbf{Quantization of Lagrangian morphisms.} For a morphism between derived algebraic stacks, 
the correct analog of a shifted symplectic structure is that of a Lagrangian structure (see 
\cite{ptvv}).  These are the maps that are candidates to 
survive after the deformation quantization. For this a version 
of the formality conjecture \ref{cform} must be stated and proved (if at all true). The basic idea here is that 
a Lagrangian morphim $f : X \longrightarrow Y$, with $Y$ $n$-shifted symplectic ($n > 0$), should deform 
$L(X)$ as an $E_{n-1}$-monoidal dg-category enriched over the deformation quantization of $L(Y)$. 
According to \cite{ca}, fully extended TQFT should be obtained this way, by quantization of fully extended TQFT with
values in a certain category of $n$-shifted symplectic derived algebraic stacks and Lagrangian correspondences
between them.

\textbf{Quantization for $n=-1,-2$.} When $n=-1$, and $n=-2$ the output of our quantization
is respectively a monoidal dg-category and braided monoidal dg-category. There are other possible
interpretations 
of the quantization in these two specific cases, as the expression "$E_{-1}$-monoidal dg-category" can 
be understood as "an object in a dg-category", and "$E_{-2}$-monoidal dg-category" as "an endomorphism of
an object in a dg-category". In particular, the quantization of a derived algebraic stack $X$ endowed
with a $(-1)$-shifted (resp. $(-2)$-shifted) Poisson structure could also be interpreted as the construction of
a deformation of an object in $L(X)$ (resp. the deformation of an endomorphism in $L(X)$). For $n=-1$ this
is the point of view taken by Joyce and his coauthors (see \cite{bbbj,bbdjs,bu}). Note that in this
setting the existence of quantization 
is predicated on the existence of orientation data which may not exist. The precise relations
with the quantization of \ref{dquant} remains to be investigated, and at the moment there is no
precise explanations of the construction of the constructible sheaf of \cite{bbbj,bbdjs,bu} 
in term of derived deformation 
theory. 

\textbf{Motivic aspects.} Deformation quantization possesses an interesing interaction with 
the motivic world. This is particularly clear when $n=-1$ (e.g. in the setting of Donaldson-Thomas theory): 
DT are made "motivic" in \cite{ks}, and the constructible sheaf $\mathcal{E}$ we mentioned above
is expected to be the Betti realization of a certain "motive" over $Bun_{G}(X)$. Because of deformation
quantization these motives most probably are instances of "non-commutative motives" over non-commutative
schemes ("$E_2$-schemes" in the setting of DT theory). For commutative base schemes non-commutative
motives have been studied in \cite{ro}, for which the constructions of \cite{bl}
provides a possible Betti realization functor. From a general point of view, the specific example
of Donaldson-Thomas theory suggests the notion of \emph{$E_n$-motives}, related to 
our deformation quantization for arbitrary values of $n$, as well as $E_n$-motives over a base
$E_{n-1}$-scheme (or stack), which is worth studying along the same lines as \cite{ro,ta2}

\textbf{Geometric quantization.} Only deformation quantization has been considered in this text. However,
derived algebraic geometry can also interact nicely with geometric quantization, a direction
currently investigated in \cite{wa}.

\end{document}